\documentclass[12pt]{amsart}
\usepackage{amssymb}
\begin{document}

\title{Alignment Correspondences}
\author{Heather Russell}

\maketitle
\tableofcontents
\date{\today}
\def\Q{{\mathbb Q}}
\def\A{{\mathbb A}}
\def\P{{\mathbb P}}
\def\C{{\mathbb C}}
\def\G{{\mathbb G}}
\def\F{{\mathbb F}}
\def\Z{{\mathbb Z}}
\def\Sp{{\textup{Spec}}}
\def\O{{\mathcal O}}
\def\H{{\textup{Hilb}}}
\def\N{{\mathbb N}} 
\def\I{{\mathcal I}}
\def\m{{\mathfrak m}}
\def\Se{{\mathcal S}}
\def\aut{{\rm Aut}}

\newtheorem{example}{Example}[section]
\newtheorem{theorem}{Theorem}[section]   
\newtheorem{corollary}{Corollary}[section]
\newtheorem{question}{Question}[section]
\newtheorem{lemma}{Lemma}[section]
\newtheorem{proposition}{Proposition}[section]

\section{Introduction}  

Let $X$ be 
a smooth, proper variety of dimension $n$ over a field $K$.  Let $R$
be the ring $K[[x_1, \dots ,x_n]]$ and $\m$ its maximal ideal.    
Let $I$
be an ideal of colength $d$ in $R$.  
The space $U(I)$ of subschemes of $X$
isomorphic to $\Sp (R/I)$ has a natural embedding in the punctual Hilbert 
scheme $\H ^d(X)$ as a locally closed subvariety (Theorem~\ref{loc}).  
More generally, one can consider the space 
$$U(I_1, \dots, I_r) = 
\{ (a_1, \dots ,a_r) \in U(I_1) \times \dots \times U(I_r):  
\exists p \in X ~{\rm and} 
$$ $$ \varphi : R \tilde{\rightarrow} {\hat \O}_{X,p} ~{\rm with}~ 
\varphi (I_1, \dots , I_r) = (a_1, \dots ,a_r)\}$$ 
and its closure $C(I_1, \dots, I_r)$ in the appropriate product of Hilbert 
schemes.  If the ideals $I_j$ are monomial, we will say that the space 
$C(I_1, \dots,I_r)$ is an {\it alignment correspondences} with 
{\it interior} 
$U(I_1, \dots , I_r)$.  The significance of the $I_j$'s being monomial
is that in this case we can show that the space $U(I_1, \dots , I_r)$
in most cases is an affine bundle on the flag bundle on $X$ and in the
remaining cases has an \'{e}tale covering by such a space and 
give a classification of these spaces via measuring sequences 
(Theorem~\ref{main}).  While many interiors of alignment correspondences are
naturally isomorphic while the alignment correspondences themselves
may vary quite a bit. In
particular, we show that in some there is a compactification
of an isomorphism class of interiors of alignment correspondences
dominating all alignment correspondences with that interior and
sometimes there is not.  Among the reasons
for studying these spaces are for their applications to enumerative
geometry, to better understand Hilbert schemes, and because they are a
fresh ground for finding something new.    

Alignment correspondences and their interiors 
have appeared in several places in the mathematical
literature previously.  If $X$ is a surface, 
given a Hilbert series $T$ space $\bar Z_T$
will be isomorphic to the fiber of $C(I)$ over $X$ 
for some ideal $I$ which may or may not be
monomial \cite{I} .  In the case that $T$ is a series of $m$
$1$'s, the ideal $I$ can be taken to be the ideal $(x_1,
x_2^m)$ and the space $C(I)$ parametrizes all zero-dimensional
subschemes of $X$ supported at a point \cite{Br}.  
The spaces $H(D)$ corresponding to an Enriques diagram $D$ 
\cite{K-P}, 
are of the form $U(I)$, although $I$ is not always a monomial ideal.
The space $U((x_1, \dots , x_{n-1}, x_n^{m+1}))$ is the space of 
curvilinear $m$-jets.  These spaces were compactified by Semple \cite{S}.  
These 
compactifications were described in modern language and coined the 
Semple bundles by Collino \cite{Col}.  Some Semple bundles can be 
realized as alignment correspondences, but it is not certain whether all 
can be.  The alignment 
correspondences $X_m = C((x,y), (x,y^2), \dots , (x,y^m))$ were
studied in \cite{Col}.  These are also 
compactifications of the spaces of curvilinear jets, but do not coincide 
with the Semple bundles for $m\ge 4$ \cite{Col}.  
Other familiar examples include flag bundles on the tangent bundle of 
$X$ and the space the jet bundles $U(x_1, \m ^n)$ parametrizing 
$m-1$-jets in the usual sense.

\vspace{.1 in}

\noindent {\bf Acknowledgements:}  I would like to thank Karen
Chandler, Joe Harris, Tony Iarrobino, Ezra Miller, Keith Pardue,
Dipendra Prasad, 
Mike Roth, Jason Starr, Ravi Vakil, and Yoachim Yameogo for their generous
help.  I would particularly like to thank Steve Kleiman for many helpful 
comments and suggestions.    

\section{Preliminaries}

\vspace{.1 in}

Let $G(I_1, \dots ,I_r)$ be the group of automorphisms of $R$ 
fixing the ideals $I_j$ for $j$ from $1$ to $r$.  The fiber of 
$U(I_1, \dots ,I_r)$ over $X$ can be identified with the quotient 
${\rm Aut}(R)/G(I_1, \dots, I_r)$ via any isomorphism 
$$\varphi: R \to \hat \O_{X,p}.$$  
If $G(I_1, \dots ,I_r) \subset G(J_1, \dots ,J_s)$, this identification
induces a 
map from $U(I_1, \dots ,I_r)$ to $U(J_1, \dots ,J_s)$.  This map 
may or may not extend to the boundary.  We will say that such maps as
well as their extensions and restrictions are {\it natural}.  
In partial compensation for
the fact that these maps do not always extend, we can consider the space 
$C( I_1, \dots ,I_r,J_1, \dots ,J_s)$ which maps to 
$C(I_1, \dots,I_r)$ and $C(J_1, \dots ,J_s)$ by projection.    
We will say this  
the space obtained by {\it superimposing} 
$C(I_1, \dots,I_r)$ and $C(J_1, \dots ,J_s)$.    
Through the isomorphism $\varphi$, the group $\aut (R)$ acts on the 
fiber of $C(I_1, \dots ,I_r)$ over $X$.  More generally if 
$C(I_1, \dots ,I_r)$ admits a natural map to $C(J_1, \dots ,J_s)$, 
assuming without loss of generality that $\varphi$ takes a point $P
\in U(J_1, \dots ,J_s)$ to $(J_1, \dots ,J_s)$  
then $G(I_1, \dots , I_s)$ acts on the fiber in $C(I_1, \dots ,I_r)$
over $P$ through $\varphi$.   

\begin{theorem}\label{loc} 
Given a sequence of ideals $I_1, \dots ,I_r$ of 
finite colengths $d_1, \dots , d_r$ respectively the space 
$U(I_1, \dots ,I_r)$, as defined in the introduction, is a locally
closed subset of the space $H= \H ^{d_1} (X)\times \dots \times \H^{d_r}(H)$.  
\end{theorem}

\noindent{\bf Proof:}  From the product of Hilbert Chow morphisms 
$$\varphi_1 \times \dots \times \varphi_r : H \rightarrow 
{\rm Sym} ^{d_1} (X)\times \dots \times {\rm Sym}^{d_r}(X),$$
we see that the points of $H$ corresponding to sequences of schemes 
all supported at the same point is a closed subvariety $Y$.  The space 
$Y$ is a fiber bundle over $X$ containing $U(I_1, \dots ,I_r)$.   
A fiber of 
$U(I_1, \dots , I_r)$ over $X$ is an $\aut (R)$ orbit of $Y$.  Hence, covering 
$X$ by neighborhoods in which there is a continuous section of 
$U(I_1, \dots , I_r)$ and a continuous action of $\aut(R)$ restricting
on the fibers over $X$ to the action described above, we see that 
the space $U(I_1, \dots , I_r)$ is locally closed in $Y$ and therefore in $H$.
$\square$

\section{Measuring Sequences}

\noindent{\bf Definition:} Let $I_1, \dots , I_r$ be a sequence of 
monomial ideals of finite colength in $R$.  
For each integer $i$ between $1$ and $n$, let 
$A_i$ be the ideal generated by images of $x_i$ under automorphisms of $R$ 
fixing each $x_j$ for $i\ne j$ and sending each $I_k$ to itself.  
We will say that the sequence 
$A_1,\dots , A_n$ is the 
{\it measuring sequence} of the sequence 
of ideals $I_1,\dots ,I_r$. 

\vspace{.1 in}

For the remainder of this paper we will continue to use the notation
of the above definition.

\begin{proposition}\label{meas0}
Each ideal $A_i$ in measuring sequence as just defined is a monomial ideal.
\end{proposition}

\noindent {\bf Proof:}
Since the 
$I_j$'s are monomial ideals, the group of 
automorphisms of $R$ fixing $x_j$ for $i\ne j$ and sending the $I_j$'s 
to themselves is 
stable under conjugation by the automorphisms scaling the $x_j$'s.  This 
together with the fact that the ideal $(x_i)$ is contained in $A_i$ make   
$A_i$ fixed by automorphisms scaling the $x_j$'s.  
Hence $A_i$ is a monomial ideal.  $\square$  

\vspace{.1 in}

The measuring sequence of a sequence of ideals is often easy to calculate.  
If the characteristic of $K$ is $0$ then we have 
$$A_i=\{f \in R : (I_j: x_i) \subset (I_j:f) ~\forall ~ 0\le j \le
r\}.$$  In a positive characteristic $p$, one must also take into
consideration what we will call the exponent types of the monomials
generating the ideals.    

\vspace{.1 in}

\noindent {\bf Definition:} If the characteristic of $K$ is a positive integer 
$p$, say that a monomial $f$ has 
{\it exponent type} $(a_0,\dots,a_m)$ if $f$ can be written as the product 
$f_0 f_1^p \dots f_m^{p^m}$ where the $f_i$'s are $p^{\rm th}$ power free 
monomials of weight $a_i$ if the weights have been given to the
variables or degree $a_i$ if no weights have been specified.  Otherwise,     
say that the exponent type of $f$ is $(wt(f))$ or, if no weights 
have been specified, $(\deg(f))$.

\vspace{.1 in} 

 
\begin{proposition}\label{meas1} 
Given a vector $$\alpha = (\alpha _1, \dots ,\alpha_n)$$ 
whose coordinates are non-negative integers
write $$\alpha =\sum _{i=0}^m p^i v_i,$$ 
where each $v_i$ is a vector with minimal non-negative integral 
coordinates and the $m$ is sufficiently large.  If the characteristic of 
$K$ is positive, let $F$ be the Frobenius map.  Let 
$$A^{\alpha} = A_1^{\alpha_1}\dots A_n^{\alpha_n}.$$  Then 
for each   
$$x^{\alpha} = x_1^{\alpha _1}\dots x_n^{\alpha _n}$$
in $I_i$, we have the containment 
$$A(\alpha) = A^{v_0}F(A^{v_1}) \dots F^m(A^{v_m}) \subset I_i.$$ 
\end{proposition}  

\noindent{\bf Proof:}  Suppose, by way of contradiction, 
that there is an exponent vector $\alpha$ and  a monomial $f$, such that 
$x^{\alpha} \in I_k$, $f \in A (\alpha) $, and  
$f \notin I_k$.  Then $f$ can be expressed as the product 
$f_0 \dots f_m^{p^m}$, where each $f_i$ is a monomial in $A^{v_i}$.  
Choose $\alpha$ so that each $f_i = x^{w_i}f_i^{\prime}$ with the 
degree of $x^{w_i}$ maximal and $f_i^{\prime} \in A^{v_i - w_i}$.  Let 
$a$ be the largest integer such that $v_i \ne w_i$.  Let $j$ be an 
integer such that the $i^{\rm th}$ coordinate $b$ of $v_a -w_a \ne 0$.  
Then there are monomials $h_1, \dots , h_b \in A_j$ such that 
$f_i^{\prime}= h_1 \dots h_bh$ with 
$h\in A^{v_a -w_a -be_j}$.  Let $g$ be the automorphism of $R$ 
with $g(x_j)= x_j + h_1 $ and $g(x_i) = x_i$ for 
$i \ne j$.   Then by Proposition~\ref{meas0}, 
$g \in G(I_1, \dots ,I_r)$.  The monomial 
$x^{\alpha}(h_1/x_j)^{p^b}$ has non-zero coefficient in $g(x^{\alpha})$ and 
hence is in $I_k$.  Let $\beta$ be the exponent vector corresponding to this 
monomial.  Then $f$ is in $A^{\beta}$, but this contradicts the maximality of 
the degrees of the $x^{w_i}$'s.  $\square$        

\vspace {.1 in}

\begin{corollary}\label{maps}  The natural map 
$$\varphi: U(A_1, \dots ,A_n) \rightarrow U(I_j)$$ is
given by sending a point $(a_1, \dots, a_n)$ to the sum ideals of the
form $a(\alpha)$ for a a set of $x^{\alpha}$'s generating $I_j$ and 
$a(\alpha)$ given by replacing each $A_i$ by $a_i$ in $A(\alpha)$ as
in Proposition~\ref{meas1}.  
\end{corollary}

\noindent{\bf Proof:} The fact that the $x^{\alpha}$'s generate $I_j$
ensures that the sum of the ideals of the form $a(\alpha)$ contains
$\varphi (I_j)$.  The fact that the sum of these ideals is contained
in $\varphi (I_j)$ follows from Proposition~\ref{meas1}. $\square$

\begin{corollary}\label{precoord} 
The group $G(A_1, \dots , A_n)$ consists of all the 
automorphisms of $R$ sending each $x_i$ to an element of $A_i$.   
\end{corollary} 

\noindent{\bf Proof:} The group $G(A_1, \dots , A_n)$ is contained in 
the set of automorphisms of $R$ sending each $x_i$ to an element of $A_i$ 
because $x_i \in A_i$.  By Proposition~\ref{meas1}, it 
follows that all such automorphisms of $R$ are contained in 
$G(A_1, \dots ,A_n)$. $\square$ 

\begin{lemma}\label{ord}  There is a bijection between nested sequences 
$$\m^2 \subsetneq B_1 \subsetneq \dots  \subsetneq B_m =\m $$ 
of distinct 
monomial ideals such that 
$G(A_1, \dots , A_n) \subset G(B_1, \dots , B_m)$, with $m$ 
is maximal and completions of the partial orderings on the
$x_i$'s such that $x_i \le x_j$ exactly when $x_i \in A_j$ to total
orderings in which no new equivalences of variables are introduced.    
\end{lemma}
\noindent{\bf Proof:} Given the set of $B_i$'s, 
let $<$ be the (possibly non-strict) 
total ordering 
on the $x_i$'s such that $x_i \le x_j$ if $x_j \in B_k$ implies 
$x_i \in B_k$.  Conversely, given such a total ordering, there is a
unique set of $B_i$'s from which it comes. The fact that the total 
ordering contains the partial ordering coming from the $A_i$'s ensures
$G(A_1, \dots , A_n) \subset G(B_1, \dots , B_m)$.  The fact that the 
$<$ is completed to a total ordering makes $m$ maximal.  $\square$

\vspace{.1 in}

Henceforth we will let $B_1, \dots B_m$ denote a sequence of ideals as
in Lemma~\ref{ord} and $B$ denote the space 
$C(B_1, \dots, B_m)$. Then $B$ is a flag 
bundle on the tangent bundle of $X$ that is maximal with the property 
that $U(A_1, \dots ,A_n)$ admits a natural map to it.  Moreover, 
we will let $G(B)$ denote $G(B_1, \dots, B_m)$, $C(I_1, \dots ,I_r,B)$
denote $C(I_1, \dots ,I_r, B_1, \dots ,B_m)$, etc.  

\begin{theorem}\label{main}  
\begin{enumerate}  
\item The space $U(A_1, \dots, A_n)$ can be expressed as 
an affine bundle over $B$.

\item The natural map $$U(A_1, \dots , A_n) \rightarrow U(I_1, \dots I_r)$$
is an \'{e}tale covering.

\item The space $C(I_1, \dots ,I_r,B)$ is 
a compactification of $U(A_1, \dots , A_n)$ which is 
a fiber bundle over $B$.    

\item $$\sum _{i=1}^n {\rm col}(A_i)= {\rm dim} (C(I_1, \dots, I_r)).$$

\end{enumerate}

\end{theorem}

\noindent{\bf Proof of Theorem~\ref{main}(1)}  
Let $W_1\oplus \dots \oplus W_m$ be a decomposition of the 
vector space $\m/ \m^2$ into direct summands spanned by monomials 
such that $B_i = W_1 \oplus \dots \oplus W_i.$  Let $N$ be large enough 
so that each $A_i$ contains $\m^N$.  For $i$ from $1$ to $m$, let 
$J_i$ be the ideal generated over $\m^N$ by the monomials spanning 
$W_i$.  Then $G(J_1, \dots ,J_m)$ is contained in $G(A_1, \dots ,A_n)$ and 
is expressible as the direct sum of the group of 
the degree preserving automorphisms 
sending the $W_i$'s to themselves with the group of automorphisms acting 
trivially on $\m/\m^N$.    By Lemma~\ref{precoord},  a set of coset 
representatives for $G(A_1, \dots ,A_n)/G(J_1, \dots ,J_m)$ is given by the 
set of all automorphisms of $R$ sending each $x_i$ to the sum 
of $x_i$, a linear combination of monomials in $A_i$ that are 
either of degree strictly between $1$ and $N$ and $x_j$'s in $A_i$ 
that are inequivalent to $x_i$.  A set of coset representatives for 
$G(B)/G(J_1, \dots ,J_m)$ is given by the set of automorphisms 
of $R$ sending each $x_i$ to $x_i$ plus a linear combination of 
monomials of degree strictly between $1$ and $n$ and $x_j$'s strictly 
less than $x_i$.  It follows that $U(J_1, \dots ,J_m)$ 
is an affine bundles over both $U(A_1, \dots ,A_n)$ and $B$.  Therefore, 
$U(A_1, \dots ,A_n)$ is an affine bundle over $B$ as well. $\square$  

\vspace{.1 in}

\begin{lemma}\label{levi} Let $W_1\oplus 0\dots \oplus W_m$ be a 
decomposition of the 
vector space $\m/ \m^2$ into nontrivial direct summands with $m$ maximal and 
each subspace $W_1 \oplus \dots \oplus W_i$ 
fixed by some power of any 
element of $G(I_1,\dots,I_r)$.  Then the group $G$ of 
degree preserving automorphisms 
of $R$ fixing the $W_i$'s is contained in $G(I_1, \dots ,I_r)$.  
\end{lemma}

\noindent{\bf Proof:} We can assume that $K$ is a perfect field 
containing an element that is neither a root of unity nor $0$, because 
if the statement of the lemma 
holds over a given field, then it also holds 
over any subfield.  By this assumption, the 
$W_i$'s are uniquely determined up to permutation by maximality of $m$ 
and spanned by monomials since they are invariant 
under scaling by some unit of $K$ that is not a root of unity.  
  
Let $M$ be the set of monic monomials in $R$.  Let the {\it $G$-topology} 
on $M$ be the topology such that the closed sets are those with span 
fixed by the group $G^{\prime} = G(I_1, \dots ,I_r) \cap G$.  
Let $B(f)$ denote the closure of a 
monomial $f$ with respect to this topology.  
We will show that the spans of the closed sets are fixed by $G$ by showing 
that the $G$-topology is the same as the exponent 
topology defined as follows.             

Let $<$ be the 
partial ordering on exponent types such that 

$$(a_0,\dots ,a_m)\le (m_0, \dots , m_b)$$
exactly when 

$$\sum_{i=0}^{k} a_ip^i \le \sum_{i=0}^k m_ip^i$$
for all $k$, with eventual equality, 
taking $a_l$ (respectively $m_l$) to be $0$ if $l>m$ 
(respectively $l>b$).  Thus only monomials of the same weight or degree 
are comparable.  

Say that a 
monomial is an {\it $i$-monomial} if it is a product of monomials in $W_i$.  
Let the {\it exponent topology} on $M$ be the 
topology on $M$ such that a set is closed if any only if whenever it contains 
a given $i$-monomial, then it contains all other $i$-monomials of lesser or 
equal exponent type.    

Let $M_i(d)$ be the set of 
$i$-monomials of degree $d$.  Expressing $M$ as the 
product $$\prod_{i,d}M_i(d)$$ both topologies are 
box topologies.  In the case of the $G$-topology, this follows from the fact 
that there is no cancellation of terms when one multiplies polynomials in  
the spans of distinct $M_i(d)$'s.  Thus we have reduced the 
problem to showing that on each $M_i(d)$ the $G$-topology and 
exponent topology are the same.  

In one direction, it is clear that the $G$-topology is as fine 
as the $p$-topology.  Suppose, by way of contradiction, that there are an  
$l$ and a $d$ such that the two topologies on $M_l(d)$ are different, 
with $d$ chosen to be minimal.  Let $B(a_1, \dots , a_t)$ 
denote the set of $l$-monomials of 
exponent type $(a_1, \dots , a_t)$ and smaller types.  
If $(a_0, \dots ,a_t)$ is not a possible exponent type, 
we will use the convention that $B(a_0,\dots, a_t)$ 
is empty.  Then there is an $l$-monomial $f$ 
of exponent type $(a_0, \dots, a_t)$ such that $B(f)$ is properly contained in 
$B(a_0, \dots ,a_t)$.  Choose $f$ to have minimal exponent type.  
Note that we cannot have  
$a_0 =0$ because then the $p^{\rm th}$ root of $f$ would give a counterexample 
of smaller degree.      
Relabelling the $x_k$'s if necessary, assume further that $f$ can be expressed 
by $$f=x_1^{e_1}\dots x_b^{e_b}h^p$$ with the degree of $h$ and then 
the $e_k$'s maximal, 
in ascending order.  Let $\alpha$ be the sequence   
$$a_0-p, a_1 -p+1 , \dots  ,a_{c-1} -p +1,a_c+1, a_{c+1}, \dots $$
where $c$ is the smallest integer with  $a_c \ne n(p-1)$.  
Then $B(\alpha)$ is the set of monomials in $M_l(d)$ of exponent type 
smaller than $(a_0,\dots ,a_t)$ not beginning with $a_0$.    
Moreover, assume that $f$ was chosen to give a vector space 
$Z \subset M_l(1)$ of maximal dimension 
such that we have the containment 
$$x_1^{e_1}\dots x_{b-1}^{e_{b-1}}Z^{e_b}h^p \subset B(f) \cup B(\alpha).$$  
  
Given $g \in G$ let $i_k$ and 
$x_k^{\prime}$ be such that $x_k^{\prime}$ differs from $g(x_k)$ by a linear 
combination of the $x_l^{\prime}$'s for $l<k$ and such that the coefficient 
of $x_{i_l}$ in $x_k^{\prime}$ is $0$ if $l<k$ and non-zero if $l=k$.    
 
By maximality of $e_1$, expanding $g(x_1^{e_1} \dots 
x_{b}^{e_b})$ out in 
terms of the $x_k^{\prime}$'s we see that modulo $B(a_0 -p, 1)$ 
the largest power 
of $x_1^{\prime}$ occurring is $e_1$ and hence that modulo 
$B(a_0 -p, 1)$, we can express 
$g(x_1^{e_1} \dots x_b^{e_b})$ as the product of $(x_1^{\prime})^{e_1}$ 
and a polynomial in the $x_k^{\prime}$'s for $k>2$.  Continuing in this 
way, we see that we have 
$$g(x_1^{e_1} \dots x_{b-1}^{e_{b-1}}) \equiv
(x_1^{\prime})^{e_1}\dots (x_{b}^\prime)^{e_{b}} \pmod {B(a_0-p,1)}.$$
Since $g$ is invertible, there must be an 
$h^{\prime}$ of exponent type $(a_1,\dots ,a_t)$ 
such that the coefficient of $h^{\prime}$ in $g(h)$ is non-zero.       
Any monomial with a non-zero coefficient in 
$$x_{i_1}^{e_1}\dots x_{i_{b-1}}^{e_{b-1}}g(Z)^{e_b} (h^{\prime})^p $$ is in 
$B(f) \cup B(\alpha)$.  This is because modulo $B(\alpha)$ there is no 
cancellation of terms when we multiply 
$g(x_1^{e_1}\dots x_{b-1}^{e_{b-1}}Z^{e_b})$ 
with $g(h^p)$.  By maximality of the dimension of $Z$, the monomial span of 
$g(Z)$ must be of the same dimension as $Z$.  Therefore, $Z$ is fixed by some 
power of any element of $G$ since there are only finitely many vector spaces 
generated by elements of $M_l(1)$ with the same dimension as $Z$.  
Hence by minimality of $W_l$ we
have $$Z=W_l.$$  By minimality of $d$, we have 
$$B(x_1^{e_1}\dots x_{b-1}^{e_{b-1}} h^p) = 
B(a_0-e_b,a_1,\dots,a_t).$$  Thus, the $G$-closure of the 
terms in  
$x_1^{e_1}\dots x_{b-1}^{e_{b-1}} W_l^b h^p$ is  
$B(a_0, \dots ,a_t)$.  Therefore we have 
$$B(f) \equiv B(a_0, \dots ,a_t) \pmod {B(\alpha)}.$$           
If $a_0\ge p $ or $p=0$, we have have reached a contradiction since then 
$B(\alpha)$ is empty.    

It remains to show that $B(\alpha)$ is contained in $B(f)$.  
By minimality of $(a_0, \dots ,a_t)$, it suffices to show that 
there is a single element of exponent type $(\alpha)$  
in $B(f)$.  Given any $l$-monomial $f_0$ of exponent type $(a_0)$, it 
follows from what we have just shown that $f_0B(a_1,\dots,a_t)^p$ is 
contained in $B(f)$.  Thus we need only show that 
$B(f_0)$ contains a monomial of exponent type $(a_0-p,1)$.  We will use 
$$f_0 =(x_1 \dots x_{i_{b-1}})^{p-1}x_b^{e_b}$$ 
where the $x_k$'s are as before, noting that for $k<b$ 
we have $e_k = p-1$.        
 
Let $Z^{\prime}$ be the maximal vector space such that the monomial span of  
$$(x_1 \dots x_{i_{b-1}})^{p-1} (Z^{\prime})^{e_b}$$ is 
contained in the span of $B(f_0)
\cup B(a_0-2p,2)$.  Without loss of generality, 
assume 
that the $x_i$'s were chosen to maximize the dimension of $Z^{\prime}$.  
Suppose, that there is an element $g\in G^{\prime}$ such that the 
monomial span of $g(Z^{\prime})$ has dimension greater than that of 
$Z^{\prime}$.    
Expanding $g((x_1 \dots x_{b-1})^{p-1}(Z^{\prime})^{e_b})$ in the 
$x_k^{\prime}$'s modulo $B(a_0 -2p,2)$, if there is a term  
of exponent type $(a_0-p,1)$, then we are done.  Otherwise, we have 
$$g((x_1 \dots x_{b-1})^{p-1}(Z^{\prime})^{e_b}) \equiv 
(x_1^{\prime}\dots x_{i_{b-1}}^{\prime})^{p-1} g(Z^{\prime})^{e_b}
\pmod{B(a_0-2p,2)}.$$  
Then since 
$$(x_1^{\prime}\dots x_{i_{b-1}}^{\prime})^{p-1} g(Z^{\prime})^{e_b}$$ 
is also contained in $B(f_0)\cup B(a_0-2p,2)$, by maximality of the 
dimension of $Z^{\prime}$, any element 
of $G^{\prime}$ preserves the 
dimension of the monomial span of $Z^{\prime}$ and hence a power 
of any element of $G$ fixes $Z^{\prime}$.  Therefore we see  
$Z^{\prime} = W_l$. 

Since there is a term of 
$(x_1^{\prime}\dots x_{i_{b-1}}^{\prime})^{p-1}(W_l)^{e_b}$ of 
exponent type $(a_0-p,1)$, 
there is a term of exponent type $(a_0-p,1)$ in $B(f_0)$.          
Hence we have arrived at a 
contradiction.  It follows that the exponent topology and the $G$-topology are 
the same.  Therefore since the action of any element of $G$ is continuous 
with respect to the exponent topology, it is continuous with respect to 
the $G$-topology and hence $G$ must be contained in 
$G(I_1, \dots ,I_r).$     $\square$ 

\vspace {.1 in}

\noindent{\bf Proof of Theorem~\ref{main}(2)}:  Throughout this proof 
we use the notation and intermediate steps from Lemma~\ref{levi}.  Let 
$B_i$ be the ideal in $R$ consisting of representatives of elements 
of $W_1 \oplus \dots \oplus W_i$.    
Assuming again that $K$ contains a unit of infinite order, 
the $W_i$'s are uniquely determined up to permutation.  Thus    
any element $g \in G(I_1,\dots,I_r)$,  
is expressible as a product $g_1g_2$ where $g_1$ preserves 
the ideals $B_i$ and $g_2$ is a permutation of the $x_i$'s sending each 
$B_i$ to $g(B_i)$.  

We will first show that $g_2 \in G(I_1, \dots ,I_r)$.  Suppose that there is a 
monomial $f$ of smallest degree such that $g_1(f)$ has $0$ coefficient 
in every element of $Gg(f)$.  Write $f = f_1f_2$ where $f_1$ is an 
$i-monomial$ for some $i$ and $f_2$ is a product of $j$-monomials for 
$j<i$.  Let $\sigma$ be the corresponding element of $S_n$.  
Then the terms of elements of 
$Gg(f_2)$ are products of $W_{\sigma (j)}$ monomials 
for $j<i$.  Thus given $g_3 \in Gg$, 
there is no cancellation of terms when multiplying the part of 
$g_3(f_1)$ spanned by $\sigma (i)$ monomials with $g_3(f_2)$.  Hence, 
by minimality of $f$, either $f_1$ or $f_2$ must be a unit.  Choose $i$ 
so that $f_2$ is a unit.  Then since $g$ is invertible, there must be 
a non-zero term in $g(f)$ that is a $\sigma (i)$-monomial of the same 
exponent type as $f$.  Thus by the proof of Lemma~\ref{levi}, all 
$\sigma (i)$ monomials of this exponent type occur in some element of 
$Gg(f)$, including $g_2(f)$.  Thus we have arrived at a contradiction 
and it follows that $g_2 \in G(I_1, \dots ,I_r)$.

Next we will show that $g_1 \in G(A_1, \dots ,A_n)$.  
The element $g_1$ has a unique expression as the product of an element 
$g_4 \in G$ and an element $g_5$ such that for each $f \in B_i$ the 
difference $f - g_5(f) \in B_{i-1}.$  By Lemma~\ref{levi} together with 
Lemma~\ref{precoord}, it follows that $g_4\in G(A_1, \dots ,A_n)$.  
    
One can see that $g_5$ is an element of $G(A_1,\dots, A_n)$ 
as follows.  Let  
$g^{\prime}\in g_5G(A_1,\dots, A_n)$ with
$$g^{\prime}(x_k) = x_k + f_k$$ and the smallest degree terms in 
the $f_k$'s of maximal degree.  Suppose that some $f_k$ has its lowest 
piece $h_k$ in $A_k$.  Then composing $g^{\prime}$ with the element 
of $G(A_1, \dots ,A_n)$ sending $x_k$ to $x_k -h_k$ and fixing the other 
$x_i$'s we get a new automorphism contradicting the maximality of the degree 
of the lowest graded piece of $f_k$.  

Suppose that $g^{\prime}$ is not the identity map.  Choose $k$ so that 
the degree of $h_k$ is minimal and if this degree is $1$, then 
$x_k$ is minimal with respect to the 
ordering $<$ associated to the $B_i$'s as in 
the proof of Theorem~\ref{main} (1).  Since 
$h_k \notin I_k$, there 
must be a monomial $f \in I_j$ for some $j$  
such that the automorphism $g_6$
sending $x_k$ to $x_k + h_k$ and fixing the other $x_i$'s does not send $f$ 
to an element of $I_j$.  However, the terms of $g_6(f)$ other that 
$f$ are the same as those in $g^{\prime}(f)$ of the corresponding 
degree having the same power of $x_k$.  Thus each monomial in $g_6 (f)$ 
and hence $g_6 (f)$ itself must 
be in $I_j$.  Thus $g^{\prime}$ is the identity map and so $g_5 \in 
G(A_1, \dots ,A_n)$.          

Thus $G(A_1,\dots, A_n)$ is the quotient of $G(I_1,\dots ,I_r)$ 
by the group of permutations of the $x_i$'s fixing the $I_j$'s 
modulo those fixing the $A_j$'s.  It follows further that 
$U(A_1, \dots, A_n)$ is an \'{e}tale covering of $U(I_1, \dots, I_r)$ 
of degree equal to the order of that group.  $\square$      

\vspace{.1 in}
\begin{lemma}\label{coord} Let $\mathcal I$ be the set of pairs 
$(i,v) $ such that $x^v$ is in the complement of $A_i$ and 
if $|v|=1$, then $x^v< x_i$ as in the proof of 
Theorem~\ref{main} (1).  Then a set of right coset representatives for 
$G(B)$ over $G(A_1, \dots ,A_n)$ is given by the 
set $\Se$ of automorphisms $g$ 
with 
$$g(x_i)= x_i + \sum_{ (i,v)\in \mathcal I} a_{i,v}x^v$$ 
for each $i$. \end{lemma}

\noindent{\bf Proof:} Let us first show that if $g \in \Se$ and 
$h \in G(A_1, \dots ,A_n)$ are automorphisms such that  
$g \circ h \in \Se ,$  then $h$
is the identity map.  Suppose that $h$
is not the identity map.  Let $k$ be such that the 
lowest graded piece $h_k$ of $h(x_i) -x_i$ is minimal.  The ordering $<$ 
can be extended to a partial ordering on monomials with the relations 
$f_1f_3 \le f_2f_4$ if $f_1 \le f_2$ and $f_3 \le f_4$.  
Then for each monomial 
$f$, the terms of $h(f) -f$ are all either of higher degree than $f$ or 
of the same degree and strictly greater than $f$.  The strictness comes from 
the consequence of Lemma~\ref{levi} that $G$ is contained in 
$G(A_1, \dots ,A_n)$. Hence the minimal 
terms of $h_k$ remains terms of $g\circ h (x_k)$.  
Thus by Lemma~\ref{precoord}, 
$g\circ h \notin \Se$.  Thus we have arrived at a contradiction and it follows 
that $h$ is the identity map.  Moreover, by the proof of 
Theorem~\ref{main}(1), the space $\Se$ is an affine space of the 
same dimension as $G(B)/G(A_1, \dots ,A_n)$ and so must 
be a full set of coset representatives.  $\square$ 

\vspace{.1 in}

\noindent{\bf Proof of Theorem~\ref{main}(3)}:  By Lemma~\ref{coord} the 
dimension of the fiber of $U(A_1^{\prime}, \dots , A_n^{\prime})$ over 
$U(B)$ is the cardinality of $\mathcal I$ in accordance with 
the theorem.  Thus it remains 
to check that $B$ has the dimension predicted by the theorem.
If $m=1$ then the measuring sequence for the sequence 
$B_1, \dots B_m$
is a sequence 
of $n$ copies of the ideal $\m$.  Thus the sum of the colengths 
of measuring ideals and the dimension of the space $C(\m)$ are 
equal to the dimension of $X$.  Suppose that the space 
$C(B_2, \dots ,B_m)$ has the expected codimension.  The measuring 
sequences of $B_1, \dots ,B_m$ and $B_2, \dots, B_m$ differ by 
${\rm cod}(B_1)$ ideals.  The differing ideals will be $B_1$ and $B_2$ 
respectively.  Hence the difference of the sums of the colengths of 
ideals in the measuring sequences is 
$${\rm cod}(B_1)( {\rm cod}(B_1) - {\rm cod}(B_2)).$$
It remains to be shown that this is the dimension of the fibers of 
$C(B_1, \dots ,B_m)$ over $C(B_2, \dots ,B_m)$.  The points 
of $C(B_2, \dots ,B_m)$ correspond to flags 
$V_m \subset \dots \subset V_2$ in the tangent bundle of $X$ where 
$V_i$ is a subspace of a tangent space of a point of dimension equal 
to the codimension of $B_i$.  The fiber over this flag consist of flags 
containing these subspaces, but with an extra subspace $V_1$ 
containing $V_2$ of dimension equal to the codimension of $B_1$.  The 
dimension of the Schubert cell parameterizing possible $V_1$'s is 
$(\dim(V_1) -\dim(V_2))(\dim(B_1))$ since $V_1/V_2$ is spanned 
by $\dim(V_1) -\dim(V_2)$ vectors taken from the quotient of the tangent 
space by $V_2$.  $\square$     

\vspace{.1 in}
   
\noindent{\bf Proof of Theorem~\ref{main}(4)}:  
Given any point $P\in B$ over a point $p\in X$,  
one can find functions $f_0, \dots, f_n$ on $X$ defined in a neighborhood 
of $p$ such that there is an isomorphism 
$$\varphi_P:\hat{R} \rightarrow \hat{\O_{X,p}} $$ with the 
following two properties.
\begin{enumerate} 
\item There are  
constants $a_{ij}$ such that we can write 
$$\varphi_P(x_i)= \sum_{j=0}^n a_{ij}f_j$$ and 
if $x_j \prec x_i$ then $a_{ij} = \delta _{ij}$.
\item The restriction of $\varphi_P$ to  $\Lambda_2$ lies  above $P$.  
\end{enumerate}
There is an open neighborhood $U$ of $P$ such 
that for each $P^{\prime} \in U$ there is a uniquely determined 
isomorphism $\varphi_{P^{\prime}}$ 
satisfying the above two properties with $P$ replaced by
$P^{\prime}$. The isomorphisms 
$\varphi_{P^{\prime}}$, being uniquely determined, vary continuously in 
$U$.  Thus about each point of $B$ we can find a 
local trivialization of $C(I_1,\dots, I_r, B)$.  
$\square$ 

\section{Universality}

\noindent{\bf Definition}: If an alignment correspondence admits a
natural map to a 
flag bundle on the tangent bundle of $X$ 
over which its interior is an affine bundle, we 
will say that the alignment correspondence is {\it directed}.  For
example the alignment correspondence $C(I_1, \dots , I_r, B)$ is
directed.   

\vspace{.1 in}

By Theorem~\ref{main}, every alignment correspondence is dominated by a 
one that is directed.  Thus for most applications, one need only 
understand directed alignment correspondences.  
Fibers of directed alignment correspondences can be understood by 
their embedding in a Grassmanian of subspaces of an appropriate 
quotient of ideals.
 
\begin{example}\label{J4}
\end{example}     
Consider the fiber $F$ of $C((x,y^4),(x,y^2))$ over 
the space $C((x,y^2))$.  By Lemma~\ref{coord}, 
the $G((x,y^2))$ orbit of 
$(x,y^4)$ is the $\A ^2$ of ideals of the form 
$$(x+ay^2 +by^3, y^4).$$  
As a subspace of the quotient $(x,y^2)/(x,y^2)^2$, 
a basis for such an ideal is 
given by   
$$\{xy+ay^3, x+ay^2 + by^3 \}.$$
Expanding $$(xy+ay^3)\wedge (x+ay^2 + b y^3) $$
with respect to the basis 
$$\{y^3 \wedge xy, y^3 \wedge y^2, y^3 \wedge x, xy \wedge y^2, 
xy \wedge x, y^2 \wedge x \}$$ of $\wedge ^2 V_4,$ 
we get a map from $F$ to $\P^5_{x_0, \dots, x_5}$ given by  
$$(x+ay^2 +by^3, y^4) \to (-b,a^2,a,a,1,0).$$
The closure of the image is cut out by  
$$x_2 - x_3 = x_5 = x_1x_4 -x_2^2 =0.$$
Thus, $F$ is a cone over a conic.  
The boundary parameterizes the $\P ^1$ of ideals of the form 
$$(\alpha xy + \beta y^2,x^2,xy^2,y^3).$$
The cone point corresponds to the ideal $$(x^2,xy,y^3).$$
$\square$ 

\vspace{.1 in}

\begin{theorem} \label {toric} Let $x_i$ have weight $e_i$ and let $g$ be 
as in Lemma~\ref{coord}.  Considering 
the $a_{i,v}$'s as coordinates for the fiber $F$ of 
$U(A_1, \dots , A_n)$ over $B$ , give 
$a_{i,v}$ weight $e_i - v$, making each $g(x_i)$ homogenious.  Then 
there are homogeneous coordinate functions 
embedding $F$ in projective space in such a way
that the closure $\bar F$ is the fiber of $C(A_1, \dots , A_n)$ over $B$.  
Moreover, if
the weights of the coordinates on $F$ are 
independent, the normalization of $\bar F $ is
a toric variety.     
\end{theorem} 

\noindent{\bf Proof:} There
is a natural embedding of $\bar F$ in 
the product a product of $n$ Grassmanians.  In particular, let $V_i$ 
be the quotient of the union of the ideals in the $G(B)$ orbit of $A_i$ 
by the the intersection of the $G(B)$ orbit of $A_i$.  Then the projection 
of $C(A_1, \dots ,A_n)$ to $C(A_i)$ gives rise to a map from 
$\bar F$ to the Grassmanian of subspaces of $V_i$ of the appropriate 
dimension.  In particular, the point of $F$ with coordinates $a_{i,v}$ 
is sent to the point corresponding to 
$g(A_i)$ viewed as a subspace of $V_i$.  These Grassmanians can then be 
embedded in projective space via homogeneous Pl\"{u}cker coordinates.  
Since $g$ preserves homogeneity, with respect to these coordinates, the 
coordinate functions will be homogeneous.  Mapping the product of projective 
spaces into a single projective space via the Segre embedding, the coordinate 
functions for $F$ are products of homogeneous coordinate functions and 
hence homogeneous.  If the weights of the coordinate functions 
are independent, then the coordinate functions can only be monomials. 
Hence the normalization 
of $\bar F$ is a toric variety.  The action of scaling the $x_i$'s gives the 
action of an open dense torus.  $\square$  

\vspace{.1 in}

\noindent{\bf Definition:}  Say that $x_i$ is {\it equivalent} to $x_j$ 
if $x_i \in A_j$ and $x_j \in A_i$.  Say that two monomials are 
{\it equivalent} if they are of the same degree in the 
variables in each equivalence class.  Say that variables $a_{i,v}$ and
$a_{j,w}$ are {\it equivalent} if $x_i$ is equivalent to $x_j$ and $x^v$ is 
equivalent to $x^w$.  

\vspace{.1 in}

\begin{theorem}\label{nonun} Suppose that there are two inequivalent 
coordinates for the fiber $F$ of $U(A_1, \dots ,A_n)$ over $B$.   
Then there is no variety $Y$ of finite type admitting a 
dominant map to every alignment correspondence with 
measuring sequence $A_1, \dots ,A_n$ such that these maps 
commute with all restrictions
of the natural maps between alignment correspondences.  

\end{theorem} 

\noindent{\bf Proof}:  Fixing the measuring sequence $A_1, \dots ,A_n$, 
it is enough to show that the theorem holds for one choice of $B$
because if $X$ exists for one choice of $B$, then $X$ is universal for
all choices of $B$.  Let $\prec$ be the ordering associated to
$B$ as in Lemma~\ref{ord}.  Without loss
of generality, assume that the indices of the $x_k$'s are such that $\prec$ 
contains the total ordering on the $x_k$'s inherited from the 
total ordering on the indices.  If there are 
inequivalent coordinates for one ordering
$\prec$, then there are inequivalent coordinates for all 
orderings.  

We first show that for some choice of 
$B$ there are two coordinates with indices given
by one of the entries of the following table.  
Suppose that these is a pair of inequivalent coordinates 
$a_{i, e_j}$ and $a_{k, e_l}$ chosen to minimize the larger of 
$i- j$ and $k-l$.  This difference must be $1$ and the indices for the pair of
coordinates must be as 
in the first or second row of the first column of the table.  
Suppose that there are no two such inequivalent coordinates.  Then there is a 
coordinate $a_{i,e_j +e_k}$ with $i$ minimal and $j$ and $k$ maximal.
Let $\prec$ be chosen so that  $i=1$ and 
$e_j +e_k$ is either equal to $2e_n$ or $e_n +e_{n-1}$.  

In the former
case, if $n \ne 2$, one can chose $B$ so there is an inequivalent 
coordinate
corresponding to the second entry in the first column of the third,
fourth, or fifth row.  If $n=2$, then there will be a second
coordinate corresponding to the entry in the sixth row.  

In the latter case, if
characteristic of $K$ is not $2$, then one can choose $B$ so that
the indices of a pair of inequivalent coordinates are given by the 
entries in the first column of the seventh row.  In characteristic
$2$, for a suitable choice of $\prec$, there will be a pair of
inequivalent coordinates corresponding to the entries in the first 
column of one of the last three rows.  Thus, one can always choose $B$ 
so that there is a pair of inequivalent coordinates with indices given by 
one of the rows of the first column of the table.  Without loss of
generality, assume that $i<j$ in the first row of the table.  In the
last two rows $N$ is an integer larger than $m$ that is $1$ less than
a power of $2$.  

\begin{table}
\begin{tabular}{|p{.75 in}| p{2.5 in}|p{.75 in}|} 
  
\hline
{$\mathcal S$}
&{$I/J_2$}
&{$J_1/I$}
\\ \hline

\noindent $i,e_{i-1} \break j, e_{j-1}$ &
$ x_i^mx_{j-1}, \dots , x_ix_{i-1}^{m-1}x_{j-1}, 
\break \break 
x_{i-1}^mx_{j}$&
$x_{i-1}^mx_{j-1} $
\\ \hline
$i,e_{i-1} \break {i+1},e_i $&
$x_i^{m+1}, \dots ,x_{i-1}^{m-1}x_i^2,x_{i-1}^{m+1} 
\break \break 
x_{i-1}^{m}x_{i+1}  $&
$x_{i-1}^mx_i$
\\ \hline
$1,2e_n \break 1,{e_{n-1}+e_n} $&
$x_1^mx_{n-1}, \dots , x_1x_n^{2m-2}x_{n-1}, 
\break \break 
x_1x_n^{2m-1} $&
$x_n^{2m}x_{n-1} $
\\ \hline 

$1, 2e_n \break 2, 2e_n $&
$x_1^m, \dots , x_1x_n^{2m-2}, 
\break \break 
x_2x_n^{2m-2}  $&
$x_n^{2m} $
\\ \hline 
$1, 2e_n \break i, e_{i-1} $&
$x_{1}^mx_{i-1}, \dots, x_1x_{i-1}x_n^{2m-2}, 
\break \break 
x_ix_n^{2m}$&
$x_{i-1}x_n^{2m} $
  \\ \hline 
$1,2e_n \break 1,3e_n $&
$x_1^{m}, \dots ,x_1x_n^{2m-2}, 
\break \break  
x_1x_n^{2m -3} $&
$x_n^{2m} $  
\\ \hline
$1,{e_{n-1}+e_n} \break n,e_{n-1} $&
$x_1^{m}, \dots ,x_1x_{n-1}^{2m-2}, 
\break \break x_{n-1}^{2m-1}x_n$&
$x_{n-1}^{2m}$
\\ \hline 
$1, {e_{n-1}+e_n} \break i, e_{i-1} $&
$x_1^mx_{i-1} x_{n-1}^{N-1} x_n^{N-1},
\dots ,x_{i-1}x_{n-1}^{N}x_n^{N},$
\break \break $x_1x_ix_{n-1}^{N-1}x_n^{N-1}$ &
$x_{i-1}x_{n-1}^{N }x_n^{N} $
  \\ \hline 
$1,{e_{n-1}+e_n} \break 2,{e_{n-1}+e_n} $&
$x_1^m x_{n-1}^{N-m}x_n^{N-m}, \dots , x_1x_{n-1}^{N-1}x_n^{N-1}, 
\break \break  x_2x_{n-1}^{N-1}x_n^{N-1} $&
$x_{n-1}^{N}x_n^{N} $
\\ \hline 

\end{tabular}
\end{table}

Fix a row of the table.  Henceforth, we will be looking only at the 
entries in that row.  Let 
$J_2$ be the ideal generated 
by monomials $f$ with weight at least that of 
all of the monomials in the second and third columns with respect to
the weighting in which $x_k$ has weight $1$ for $x_a \prec x_k$ 
(not excluding $x_k \prec x_a$) and weight $2$ otherwise, and either 
the weight of $f$ or the weighted exponent type is strictly greater than
those of each these monomials for some $a$.  
Let $I_m$ be the ideal generated over $J_2$ by 
the monomials of the second column and equivalent monomials.
Then $J_2$ is the intersection of the ideals in the $G(B)$ orbit
of $I_m$ and these monomials are a basis for $I_m/J_2$.  Moreover, these
monomials together with the monomials equivalent to the one in the
third column are a basis for the union $J_1$ of the ideals in
the $G(B)$ orbit of $I_m$ over $J_2$.  Let $\Se$ the set of elements in 
the first column.  Let $L$ be the lattice generated by the 
vectors $e_i - v$ for $i,v \in \Se$.  Let 
$\mathcal A$ 
denote the set of automorphisms $g$ of $R$ with 
$$g(x_i) = x_i + 
\sum_{(i,\alpha)\in {\Se}}a_{i,\alpha}x^{\alpha}$$
for each $i$ from $1$ to $n$.
Let $M$ be the subspace of $J_1/J_2$ generated by the monomials in the
second and third columns.  Then the differences of the weights of
these monomials lie in $L$.  

The closure of the $\mathcal A$ 
orbit of $I_m$ in the Grassmanian of hyper-planes in $J_1/J_2$ 
is a cone over a rational normal curve of degree $m$.  
The $\mathcal A$ orbit of $(I_2, \dots , I_m)$ 
in the appropriate product of Grassmanians.  The boundary
of the image is a chain of $(m-1)$ $\P^1$'s.  

Suppose there is a
variety $Y^{\prime}$ with a dominant maps $f_m$ to the fiber   
of $C(I_2, \dots ,I_m, B)$ over $B$ for all $m$ such that these maps
commute with the natural projection maps.  If $Y$ exists, then 
its fiber over a point in $B$ is such a variety.  Then the image of the
boundary of the fiber of $C(I_2, \dots ,I_m,B)$ under $f_m$ 
is independent of $m$ and has at least $m-1$ components.  Thus $Y^{\prime}$
cannot be of finite type.  Thus $Y$ cannot be of finite type. $\square$

\begin{theorem}\label{exceptions}
If the measuring sequence 
$A_1,\dots A_n$ does not satisfy the hypotheses of Theorem~\ref{nonun} 
then the coordinates for the fiber $F$ of $U(A_1, \dots, A_n,B)$ over $B$ 
are all coordinates 
$a_{i,v}$ in an equivalence class satisfying one of the following:
\begin{enumerate} 
\item $v = e_j$,
\item $v= e_j+e_k$ with $x_j$ and $x_k$ equivalent.   
\item $v= e_j+e_k$ with $x_j$ and $x_k$ equivalent and $j \ne k$   
\item $v= e_j+e_k$ with $x_k$ in a unique equivalence class and $j \ne k$.   
\end{enumerate}  
The last two cases happen only in characteristic $2$.  
Let $m$ be the number of elements in the equivalence class of $x_i$ 
and $l$ the number in the equivalence class of $x_j$ for any
coordinate $a_{i,v}$.  Then if $m$ or
$l= 1$ in the first or fourth case, $l=1$ in the second, or $l=2$ in
the third, then given a sequence $(I_1, \dots, I_r)$
with measuring sequence 
$(A_1, \dots, A_n)$, the space $C(A_1, \dots, A_n,B)$ dominates 
$C(I_1, \dots, I_r)$ via the natural map.   
\end{theorem}

\noindent{\bf Proof:}  The fact that these are the only cases follows from
the proof of Theorem~\ref{nonun}.  Let $G$ be the group of linear
automorphism of $R$ fixing the subspace $V_1$ of $\m/ \m^2$ 
spanned by the variables
equivalent to $x_i$, fixing the subspace $V_2$ 
spanned by variables equivalent to
$x_j$, and fixing all all variables.   
In the first and last case, $F$ can
be identified with the unipotent radical $U_{m,l}$ 
of the parabolic subgroup $P_{m,l}$ of $GL(V_1 \oplus V_2)$ fixing $V_1$.  
The group $G$ can be identified with the Levi subgroup of $P_{m,l}$
acting by conjugation.  In the second case 
$F$ can be identified with the space of $m$-tuples
of quadratic forms in $l$ variables.  
An element of $G$ viewed as an element $(g_1,g_2) \in GL(V_1) \times
GL(V_2)$ acts by the composition
of the action of $g_1$ by change of 
variables and the action of $g_2$ by matrix 
multiplication.  The third case is similar to the second, except that
the space of $m$-tuples of quadratic forms is taken modulo the image of
the Frobenius on the space of linear forms.   

If $m$ and $l$ are as in the hypotheses of the theorem then $G$ is the
group of all linear transformations of the affine space $F$.  In the third case, this is less straightforward, 
but using the fact that the characteristic is $2$, 
$GL(V_2)$ acts on $F$ by the cofactor matrix of the matrix
giving the action on $V_2$.  
The fiber of $C(A_1, \dots,A_n,B)$ over $B$ is a projective space in which
the only fixed subvariety under $G$ is the boundary.  Thus it follows
that it is the universal $G$-equivariant
compactification of $F$ and hence that 
$C(A_1, \dots,A_n,B)$ is universal. $\square$

\vspace{.1 in}

In the four types of measuring sequences listed in
Theorem~\ref{exceptions}, the coset representatives for $G(B)/G(A_1,
\dots ,A_n)$ given in Lemma~\ref{coord} form an algebraic group.  Thus
equivariant compactrifications of these groups arise naturally as 
fibers of alignment correspondences with these measuring sequences over $B$.


\begin{thebibliography}{60}



\bibitem  {Br}Brian\c{c}on, J. 
{\it Description de $\H ^n \C \{x,y\}.$} Inventiones math. 
{\bf 41},45--89 (1977).




\bibitem{Col}Collino, A. {\it Evidence for a conjecture of Ellingsrud 
and Str{\o}mme on the Chow ring of $\H _d(\P^2)$.}  Illinois Journal of 
Mathematics {\bf 32}, 171--210 (1988).



\bibitem {F}Fulton, W. {\it Introduction to Toric Varieties}, 
Princeton University Press, Princeton, 1993.

\bibitem {I}Iarrobino, A. {\it Punctual Hilbert schemes} Mem. Amer. Math. Soc. 
{\bf 10} no. 188 (1977).

\bibitem {K-P} Kleiman, S., R. Piene {\it Enumerating singular curves on 
surfaces}, Algebraic Geometry: Hirzebruch 70 (Warsaw 1998), 209--238, 
Contemp. Math., 241, Amer. Math. Soc, Providence, RI, 1999.

\bibitem {N}Nakajima, H. {\it Lectures on Hilbert schemes of points 
on surfaces}, University Lecture Series, {\bf 18}, Amer. Math. Soc.,
Providence, RI, 1999.

\bibitem {R1}Russell, H. {\it The Enumeration of Plane Curves 
with Singularities Corresponding to Monomial Ideals} (preprint)

\bibitem {R2}Russell, H. {\it Hilbert Schemes and Monomial Ideals},
thesis, Harvard University, June 1999.

\bibitem {S} Semple, J. { \em Some investigations in the geometry of curves 
and surface elements}. Proc. London Math. Soc. (3), vol. 4 (1954), pp 24-49.

\bibitem {V} Vakil, R. {\it A beginner's guide to jets} (preprint)

\end{thebibliography}
\end{document}